\documentclass[12pt]{article}
\usepackage{amssymb,amsfonts,amsmath}
\usepackage{graphics}
\usepackage{graphicx}
\usepackage{epsfig}
\usepackage{color}
\usepackage{enumerate}

\newcommand\CC{{\mathbb C}}

\newcommand\DD{{\mathbb D}}
\newcommand\NN{{\mathbb N}}
\newcommand\RR{{\mathbb R}}

\newcommand\TT{{\mathbb T}}

\def\beq{\begin{equation}}
\def\eeq{\end{equation}}
\newtheorem{thm}{Theorem}[section]
\newtheorem{prop}[thm]{Proposition}
\newtheorem{defn}[thm]{Definition}
\newtheorem{lem}[thm]{Lemma}
\newtheorem{cor}[thm]{Corollary}
\newtheorem{rem}[thm]{Remark}

\newtheorem{ex}[thm]{Example}
\newtheorem{Not}[thm]{Notation}

\newcommand\beginpf{\noindent {\bf Proof:} \quad}

\newcommand\re{\mathop{\rm Re}\nolimits}

\newcommand\esssup{\mathop{\rm ess\ sup}\nolimits}

\def\beginpf{\noindent {\bf Proof:} \quad}
\def\endpf{\rightline{$\square$}}

\def\N{{\mathbb N}}

\def\T{{\mathbb T}}

\newcommand\RE{\mathop{\rm Re}\nolimits}
\newcommand\IM{\mathop{\rm Im}\nolimits}
\renewcommand\phi{\varphi}

\newcommand\D{\mathcal {D}}

\title{A class of quasicontractive semigroups acting on Hardy and Dirichlet space}
\author{ C. Avicou\thanks{I. C. J., UFR de
Math\'ematiques,  Universit\'e  Lyon~1, 43 bld. du 11/11/1918,
69622 Villeurbanne Cedex, France.
  \protect\linebreak[3]
{\tt avicou@math.univ-lyon1.fr}.},
 I. Chalendar\thanks{I. C. J., UFR de
Math\'ematiques,  Universit\'e  Lyon~1, 43 bld. du 11/11/1918,
69622 Villeurbanne Cedex, France.
  \protect\linebreak[3]
{\tt chalendar@}math.univ-lyon1.fr}
\ and   J.R. Partington\thanks{School of
Mathematics,
University of Leeds, Leeds LS2 9JT, U.K.
\protect\linebreak[3]
{\tt J.R.Partington@leeds.ac.uk}.}
}

\begin{document}

\baselineskip18pt

\maketitle
\begin{abstract}
This paper provides  a complete characterization of quasicontractive $C_0$-semigroups on  
Hardy and  Dirichlet space with a prescribed generator of the form $Af=Gf'$. 
We show that such semigroups are semigroups of composition operators and 
we give simple sufficient and necessary condition on $G$.  
Our techniques are based on ideas from semigroup theory, such as 
the use of numerical ranges.
\end{abstract}

\noindent\textsc{Mathematics Subject Classification} (2000):
Primary: 47D03, 47B33
Secondary:  47B44, 30H10

\noindent\textsc{Keywords}:
strongly continuous semigroup,  Hardy space, Dirichlet space,
quasicontractive semigroup, composition operators.

\bibliographystyle{plain}
%\newpage
%\tableofcontents 

\section{Introduction}
In this paper, an operator is always assumed to be  linear but not necessarily
bounded.  

Let $X$ be a Banach space. A one-parameter family $(T(t))_{t\geq 0}$ of bounded linear operators from $X$ to $X$ is a semigroup of bounded 
linear operators on $X$ if 
\begin{enumerate}[(i)]  
\item $T(0)=\mathop{\rm Id}$,  the identity operator on $X$;
\item $T(t+s)=T(t)T(s)$ for every $s,t\geq 0$.
\end{enumerate}
The linear operator $A$ defined by 
\[D(A)=\left\{ x\in X: \lim_{t\downarrow 0}\frac{T(t)x-x}{t}\mbox{ exists}\right\}\]
and 
\[Ax=\lim_{t\downarrow 0}\frac{T(t)x-x}{t}\mbox{ for }x\in D(A)\]
is the (infinitesimal) generator of the semigroup $T(t)$, $D(A)$ is the domain of $A$. 
The semigroup $(T(t))_{t\geq 0}$ of bounded linear operators on $X$ is strongly continuous if 
\[\lim_{t\downarrow 0}\|T(t)x-x\|_X=0 \mbox{ for every }x\in X.\] 
Such semigroups are also called $C_0$-semigroups. 

A straightforward consequence of the uniform boundedness theorem is that 
given a $C_0$-semigroup $(T(t))_{t\geq 0}$ on the Banach space $X$, 
there exist $w\geq 0$ and $M\geq 1$ such that 
\[\|T(t)\|\leq Me^{wt}\mbox{ for all }0\leq t<\infty.\]
In the particular case where $M=1$, the semigroup is said to be
quasicontractive. For $M=1$ and $w=0$, $(T(t))_{t\geq 0}$ is a
semigroup of contractions.   \\

In 1978, Berkson and Porta \cite{B_P} gave a complete description of
the generator $A$ of semigroups of composition operators on the Hardy
space $H^2(\DD)$ (see Section~\ref{sec:2.3}),
induced by a semigroup of analytic self-maps of $\DD$ (see
Section~\ref{sec:2.2} for the definition of such semigroups). 
Abate \cite{Abate} rediscovered the main results of \cite{B_P}, using a
different approach, and considering higher dimensions of the scalar
space. 
 
 Berkson and Porta \cite{B_P} noticed that such semigroups are strongly continuous on $H^2(\DD)$ and 
Siskakis
\cite{S, sisk96} noticed that they are
strongly continuous on the Dirichlet space $\D$. 
Moreover, it is not difficult to see that the generator $A$ of
a semigroup of composition operators is of the form $Af=Gf'$.   The
aim of this paper is to give a complete description of quasicontractive
$C_0$-semigroups of bounded operators on $H^2(\DD)$ and $\D$ whose
generator $A$ is of the form $Af=Gf'$; that is, unlike previous authors, we do not
assume {\em a priori\/} that we are working with semigroups of composition operators..\\ 

The paper is organized as follows. First, in Section~\ref{sec:2}, we
recall the Lumer--Phillips theorem in order to obtain a contractive or
quasicontractive 
$C_0$-semigroup by means of the numerical range of its generator. We
also present the main result of \cite{B_P} concerning the semigroups
of holomorphic functions on $\DD$. Then, the weighted Hardy spaces are
defined, their main properties are recalled (in particular the fact
that some of them are reproducing kernel Hilbert spaces is emphasized) and finally we also study
the norm of composition operators induced by a univalent symbol
$\varphi$ on $H^2(\DD)$ and $\D$. This study is essential to check
that the semigroups of composition operators are indeed strongly
continous and quasicontractive. To that aim, on the Dirichlet space,
the optimal  estimates proved in \cite{martin05} are crucial.  

In Section~\ref{sec:3}, we present our main result on $H^2(\DD)$,  Theorem~\ref{ThmB},
 which asserts that the only quasicontractive
$C_0$-semigroups whose generator is of the form $Af=Gf'$ are the
semigroups of composition operators. We also prove necessary and
sufficient conditions on $G$, different from the one of Berkson and
Porta.

In Section~\ref{sec:4}, we prove in Theorem~\ref{thm:dir} that the assertions of
Theorem~\ref{ThmB} are equivalent to the fact that $A$ generates a
quasicontrative $C_0$-semigroup on the Dirichlet space, which is
itself equivalent to the fact that $A$ generates a semigroup of
composition operators on $\D$.  

The last section contains comments on well-known algorithms to test
the conditions on $G$, and explicit examples of constructions of the
semigroup of composition operators for a class of analytic polynomial $G$.

\section{General background}\label{sec:2}
\subsection{Characterization of contractive $C_0$-semigroups of
  bounded operators on a Hilbert space}
 
 Besides the well-known Hille--Yosida theorem (see for example Thm.~3.1 in \cite{Pazy})  which characterizes  $C_0$-semigroups in terms of the growth of the resolvent of their generators $A$, another  
 useful theorem  is the Lumer--Phillips theorem (see for example Thm.~4.3 in \cite{Pazy}) which is well-adapted for the characterization of quasicontractive $C_0$-semigroups in terms of the numerical range of $A$. 
 
 From now on, we assume that the Banach space $X$ on which $T(t)$ is defined is a complex Hilbert space,
and we denote it by $H$. This hypothesis will simplify the definition of dissipative operators involved in the Lumer--Phillips theorem.

 Let $A:D(A)\to H$ be a linear operator. Then $A$ is dissipative if 
 \[\RE \langle Ax,x\rangle \leq 0\mbox{ for all }x\in D(A), \|x\|=1.\]
In other words, $A:D(A)\to H$ is dissipative if  the numerical range of $A$ is in the left half-plane.  
 
 \begin{thm}[Lumer--Phillips]
 Let $A$ be a linear operator with dense domain $D(A)$ in $X$. 
\begin{enumerate}[(i)]
\item If $A$ is dissipative and there exists $\lambda_0>0$ such that
  $(\lambda_0\mathop{\rm Id}-A)D(A)=H$, then $A$ is the generator of
  contractive $C_0$-semigroup 
on $H$.   
\item If $A$ is the generator of a contractive $C_0$-semigroup on $H$,
  then $A$ is dissipative and for all $\lambda >0$, 
$(\lambda\mathop{\rm Id}-A)D(A)=H$. 
\end{enumerate}
\end{thm}
 This theorem is also of great interest for the characterization of
 quasicontractive $C_0$- semigroups observing that $\|T(t)\|\leq e^{wt}$ if and only if 
 $\|\widetilde{T}(t)\|\leq 1$, where  $\widetilde{T}(t):=T(t)e^{-wt}$ is the semigroup whose generator is $A-w\mathop{\rm Id}$, if $A$ is the generator of 
 $(T(t))_{t\geq 0}$.  
In particular we have then the following result. 
\begin{cor}\label{cor:lm}
Let $A:D(A)\to H$ be a linear operator with a dense domain. Then $A$ generates a quasicontractive
$C_0$-semigroup if and only if $\sup\{\RE\left( \langle
  Ax,x\rangle\right):x\in D(A), \|x\|=1\}<\infty$ and there exists
$\lambda>0$ such that $(A-\lambda\mathop{\rm Id})D(A)=H$. 
\end{cor}

%\section{Semigroup of composition operators on $H^2(\DD)$}
%\section{General background}
\subsection{Semigroups of analytic functions}\label{sec:2.2}
\begin{thm}[Denjoy-Wolff]
Let $\varphi:\DD\to\DD$ holomorphic such that $\varphi$ is not an elliptic automorphism. Then there is a point 
$b\in \overline{\DD}$ such that $\varphi_n$(:=$\varphi\circ\varphi\circ\cdots\varphi$, $n$ times) converges to $b$ uniformly 
on compact subsets of $\DD$. \\
If $|b|<1$, then $\varphi(b)=b$, while if $|b|=1$, then $b$ behaves as a fixed point in the sense that $\lim_{r\to 1^-}\varphi (rb)=b$. 
This distinguished point is called the Denjoy--Wolff point of $\varphi$.  
\end{thm}
In the exceptional case of elliptic automorphisms different from the identity map, the sequence of iterates move around 
an interior fixed point without converging to it. 
\begin{defn}
A one-parameter semigroup of analytic functions of $\DD$ into
itself is a family $\Phi=\{\varphi_t:t\geq 0\}$ of analytic self-maps of $\DD$
such that 
\begin{enumerate}
\item $\varphi_0(z)=z$ for all $z\in \DD$;
\item $\varphi_{t+s}(z)=\varphi_t\circ \varphi_s(z)$ for all $t,s\geq
  0$ and $z\in \DD$;
\item $(t,z)\mapsto \varphi_t(z)$ is continuous on $[0,\infty)\times
  \DD$. 
\end{enumerate}

\end{defn}
Using Vitali's Theorem on convergence of holomorphic functions, it
follows that the continuity of $(t,z)\mapsto \varphi_t(z)$  on $[0,\infty)\times
  \DD$  is equivalent to the continuity of $t \mapsto \varphi_t(z)$
for each $z\in \DD$. 
%Then one can define the semigroup $T$ of composition operators by,
%$\forall t\in\RR_+,\ T(t)=C_{\varphi_t}:f\m%apsto f\circ\varphi_t$. 
Such semigroups have  extensively been studied by
Berkson and Porta \cite{B_P} (see also \cite{AS}), who proved the 
following useful result.
 %We define the semigroup $T$ of composition operators by, $\forall t\in\RR_+,\ T(t)=C_{\varphi_t}:f\mapsto f\cir%c\varphi_t$. Such semigroups have already been studied (see \cite{B_P} or \cite{AS}).\\
%We may recall the main results needed:
\begin{prop}\label{Sisk}
Let $\Phi=(\varphi_t)_{t\in\RR_+}$ be a semigroup of analytic functions on $\DD$, then:
\begin{enumerate}[(i)]
\item For every $t\in\RR_+$, the function $\varphi_t$ is univalent.
\item There is a holomorphic mapping $G:\DD\to \CC$ called the
  generator of $\Phi$ such that 
\begin{equation}\label{eq:gene}
\frac{\partial \varphi_t(z)}{\partial t}=G(\varphi_t(z))
\end{equation}
for all $t\geq 0$ and all $z\in \DD$. The convergence  
\[G(z)=\lim_{t\to 0^+} \frac{\partial \varphi_t(z)}{\partial t}\]
is uniform on every compact subsets of $\DD$.
\item Moreover, the infinitesimal generator $G$ of $\Phi$ has the unique representation
\[G(z) = F(z)(\overline{\alpha}z-1)(z-\alpha),\ \forall z\in\DD,\]
where $F:\DD\to\CC$ is analytic and satisfies $\RE(F)\geqslant 0$, and $\alpha$ is the Denjoy--Wolff point of one (and thus any) $\varphi_t$, $t>0$.
\item Conversely, let $F:\DD\to\CC$ be analytic with $\RE(F)\geqslant 0$ and $\alpha\in\overline{\DD}$, the function 
$z \mapsto F(z)(\overline{\alpha}z-1)(z-\alpha)$, generates a semigroup of analytic function on $\DD$.
%\item[.] Moreover the infinitesimal generator $\Gamma$ of $T$ is given by
%\[\Gamma f(z) = G(z)f'(z).\]
\end{enumerate}
\end{prop}
% \textbf{KEEP or NOT}
% The previous expression of $G$ can be studied in two cases according to the Denjoy--Wolff point: when $|\alpha|<1$, up to a M\"obius transformation, we can equivalently consider that $\alpha$ is $0$, in that case, we shall say that $\Phi$ belongs to the class $\Phi_0$; when $|\alpha|=1$, up to a rotation, this reduces to $ \alpha=1$ and we will then say that $\Phi$ belongs to the class $\Phi_1$.
\subsection{Operators on weighted Hardy spaces}\label{sec:2.3}
Our aim is to study semigroups of bounded operators on classical spaces of
analytic functions such as the Hardy space $H^2(\DD)$ and the
Dirichlet space $\cal D$, which are particular cases of the so-called
``weighted Hardy spaces''.  
\begin{defn}
Take $(\beta_n)_{n \ge 0}$ a sequence of positive real numbers.
Then $H^2(\beta)$ is the space of analytic functions 
\[
f(z)=\sum_{n=0}^\infty c_n z^n
\]
 in the unit disc $\DD$ that have finite norm
\[
\|f\|_\beta= \left( \sum_{n=0}^\infty |c_n|^2 \beta_n^2 \right)^{1/2}.
\]
The case $\beta_n=1$ gives the usual \textbf{Hardy} space $H^2 (\DD)$. \\
The case $\beta_0=1$ and $\beta_n=\sqrt{n}$ for $n\geq 1$ provides the \textbf{Dirichlet} space
$\cal D$, which is included in $H^2(\DD)$. \\
The case $\beta_n=1/\sqrt{n+1}$ produces the \textbf{Bergman} space,
which contains $H^2(\DD)$.
\end{defn}

Obviously the polynomials are in $H^2(\beta)$ and with an extra condition on $(\beta_n)_n$, the Hilbert space $H^2(\beta)$ is also a 
{\em reproducing kernel space}, i.e. for all $w\in \DD$, there exists a function $k_w\in H^2(\beta)$ such that 
\[\langle f,k_w\rangle =f(w),\]
for all $f\in H^2(\beta)$ (see p. 19 in \cite{Cowen-Mac-Cluer} and p. 146 in \cite{Par97}). More precisely, if $(\beta_n)_n$ is such that 
\begin{equation}\label{eq:rks}
\sum_{n\geq 0}\frac{|w|^{2n}}{\beta_n^2} <\infty\mbox{ for all }w\in\DD,
\end{equation}
it follows that $H^2(\beta)$ is a reproducing kernel Hilbert space and 
\[k_w(z)=\sum_{n\geq 0}\frac{\overline{w}^n}{\beta_n^2}z^n \qquad \mbox{with}
\qquad \|k_w\|^2_{H^2(\beta)}=\sum_{n\geq 0}\frac{|w|^{2n}}{\beta_n^2}.\] 
In fact (\ref{eq:rks}) is also equivalent to the more explicit condition  $\liminf (\beta_n)^{1/n}\geq 1$.

Given an operator $A$ (possibly unbounded) defined by $A f(z) = G(z) f'(z)$ on its domain $D(A)=\{f\in H^2(\beta),\ Gf'\in H^2(\beta)\}$ where $G\in H^2 (\beta)$, we  would like to know if there exists a  $C_0$-semigroup on $H^2(\beta)$ with generator $A$. The next proposition asserts that two necessary conditions for $A$ to be a $C_0$-semigroup generator are satisfied.

\begin{prop}
Let $(\beta_n)_{n \ge 0}$ a sequence of positive real numbers such that $zH^2(\beta)\subset H^2(\beta)$.
Any operator $A$ defined by $A f(z) = G(z) f'(z)$ on $D(A)$ where $G\in H^2(\beta)$ is densely defined on $H^2(\beta)$ and closed.
\end{prop}

\beginpf
The operator $A$ is defined on polynomials which form a dense family of function in $H^2(\beta)$.
We consider a sequence $(f_n: z \mapsto\sum_k a_k^nz^k)\in H^2(\beta)^\N$  and two functions $f: z \mapsto\sum_k a_kz^k,
g: z \mapsto\sum_k c_kz^k \in H^2(\beta)$ such that $f_n \rightarrow f$ and $Gf_n' \rightarrow g$ in $H^2(\beta)$. We denote $G(z)=\sum_k b_kz^k$. 
We now consider the truncated sums, up to the $N$-th exponent:
\[\|(Gf_n' - Gf')_N \|_2^2 \leqslant  \sum_{k=0}^N \beta_k\left|\sum_{i+j=k} b_{i-1}j(a_{j}^n -a_j)\right|^2.\]
As  $f_n \rightarrow f$ in $H^2(\beta)$, one has 
\[ \sum_{k=0}^\infty \beta_k|a_{k}^n -a_k|^2 \rightarrow 0\]
and thus $\forall k\in\N, |a_{k}^n -a_k| \rightarrow 0$.
Hence, 
\[
\|(g-Gf')_N\|_2 \leqslant \|(g-Gf_n')_N\|_2+ \|(Gf_n' - Gf')_N \|_2\rightarrow 0.
\]
We have shown that $\forall k\leqslant N$, $c_k=\sum_{i+j=k}j
b_{i-1}a_j$ and since this can be done for each choice of $N$, we
conclude that $g=Gf'\in H^2(\beta)$, $Gf_n' \rightarrow Gf'$ in
$H^2(\beta)$ and $A$ is closed.

\endpf

\begin{prop}
If $A$ is defined by $A f(z) = G(z) f'(z)$ where $G\in H^2(\DD)$ and
$\frac 1G\in H^\infty(\DD)$, then
 $A $ cannot be the generator of a one-parameter semigroup.
\end{prop}
\beginpf
Let $\lambda$ be a real number; then $\lambda\in\sigma(A)$ if  there exists $f\in H^2(\DD)$ such that
\[ G(z)f'(z) =\lambda f(z).\]
Since $\frac 1G\in H^\infty(\DD)$ the function $u=\int \frac{\lambda}{G}dz$ lies in $H^\infty(\DD)$ and $f=e^{u}\in H^\infty \subset H^2(\DD)$ satisfies $G(z)f'(z) =\lambda f(z)$. Thus $\RR \subset \sigma(A)$. This cannot occur for $C_0$-semigroups, see e.g. \cite[Chap.II,1.13]{E_N}.

\endpf
\begin{cor}
If $A$ is defined on $H^2(\DD)$ by $A f(z) = p(z) f'(z)$ where $p$ is polynomial with no roots in the unit closed disc $\overline{\DD}$, then $A$ cannot be the generator of a one-parameter semigroup.
\end{cor}

\subsection{Bounded composition operators on Hardy  and 
  Di\-ri\-chlet spaces}
Composition operators on the Hardy space $H^2(\DD)$ have a quite surprising
property, namely, provided that they are well-defined, they are
always continuous. This fact is not true on the Dirichlet
space. Moreover, we have the following upper bound for the norm (see
Thm.~3.8 in \cite{Cowen-Mac-Cluer}).
\begin{thm}\label{th:compoH2}
Let $\varphi:\DD\to \DD$ be an analytic function. Then
$C_\varphi$ maps $H^2(\DD)$  continuously into $H^2(\DD)$,  and moreover 
\[\| C_\varphi  \|\leq \left(
  \frac{1+|\varphi(0)|}{1-|\varphi(0)|}\right)^{1/2}.\]
\end{thm}

The previous result is very useful to estimate the growth of the norm
of semigroups of composition operators  on the Hardy space. Indeed, a
first consequence of 
Theorem~\ref{th:compoH2} is that each semigroup $\Phi$  of analytic functions
on $\DD$ induces a $C_0$-semigroup of bounded operators on $H^2(\DD)$. 

\begin{cor}\label{cor:quasicont}
Let $\Phi=(\varphi_t)_{t\geq 0}$ be a semigroup of  analytic functions
on $\DD$. Then $(C_{\varphi_t})_{t\geq 0}$ is a quasicontractive $C_0$-semigroup   on $H^2(\DD)$. 
\end{cor}

\beginpf
The continuity of $t\mapsto \varphi_t(0)$ implies that $K:=\{\varphi_t(0):0\leq t\leq 1\}$ is a compact subset of 
$\DD$. Since $G$, the generator of $\Phi$, is holomorphic on $\DD$, we get 
\[\sup_{0\leq t\leq 1}|G(\varphi_t(0))|<\infty.\]
By (\ref{eq:gene}), 
it follows that there exists $M>0$ such that $|\varphi_t(0)|\leq Mt$, and then, 
for $0\leq t\leq \frac{1}{2M}$, $|\varphi_t(0)|\leq \frac{1}{2}$. 
Using Theorem~\ref{th:compoH2}, we also know that 
\[ \|C_{\varphi_t}\|\leq \left( \frac{1+|\varphi_t(0)|}{1-|\varphi_t(0)|}\right)^{1/2},\]
which implies that $\|C_{\varphi_t}\|\leq 1+O(t)$, and thus there exists $w\geq 0$ such that 
  \[\| C_{\varphi_t}  \|\leq e^{wt},\]
for all $t\geq 0$. 
Moreover the hypotheses on $\Phi$ imply that $C_{\varphi_t}f(z)$ tends to $f(z)$ as $t$ tends to $0$, for all $z\in\DD$ and all 
$f\in H^2(\DD)$. In other words the semigroup $(C_{\varphi_t})_t$ is weakly continuous. It follows that 
$(C_{\varphi_t})_t$ is strongly continuous (see Thm. I.5.8 in \cite{E_N}) .

\endpf

On the Dirichlet space, it is not true that $C_\phi$ is well defined 
whenever $\varphi$ is a self-map of $\DD$. For example, for $\phi$ an infinite Blaschke
product, $C_\varphi$ is not a bounded composition operator on $\cal
D$. Nevertheless, if $\varphi$ is univalent, $C_\varphi$ is bounded on
$\cal D$ (see Section~6.2 of \cite{Ransfordetal}). We have therefore
the following preliminary result.
\begin{prop}\label{prop:compodirichlet}
Let $\Phi=(\varphi_t)_{t\geq 0}$ be a semigroup of  analytic functions
on $\DD$. Then $(C_{\varphi_t})_{t\geq 0}$ is a semigroup of bounded
operators on $\cal D$. 
\end{prop}
Using \cite{CP14}, it is possible to prove that
Proposition~\ref{prop:compodirichlet} is still true for any space
$H^2(\beta)$ containing $\cal D$. But this is beyond the scope of this
paper.

\section{Quasicontractive semigroups on the Hardy space}\label{sec:3}

From now on, the function $G$ will lie in $H^2(\DD)$ and the operator $A$ will be defined by $A f=G f'$ on the domain $D(A)=\{f\in H^2(\DD),\ Gf'\in H^2(\DD)\}$.

\begin{thm}\label{ThmB}
The operator $A$ generates a $C_0$-semigroup of composition operators on $H^2(\DD)$ if and only if $\forall z\in\DD$,
\begin{equation}\label{B}
2\RE(\overline{z}G(z)) + \left( 1 - |z|^2\right)\RE(G'(z))\leqslant 0.
\end{equation}
\end{thm}
\beginpf
 Suppose $A$ is such a generator, let $(\varphi_t)$ denote the corresponding semigroup. From analyticity, one has for small $t$ and fixed $z$:
\[ \varphi_t(z) = z+G(z)t+o(t),\]
\[ \varphi'_t(z)= 1 +G'(z)t+o(t).\]
From the Schwarz--Pick lemma (see \cite{A}),
\[ |\varphi_t'(z)| \leqslant \frac{ 1-|\varphi_t(z)|^2}{1-|z|^2},\]
and thus,
\[1+\RE(G'(z))t+o(t)\leqslant \frac{1-|z|^2}{1-|z|^2}-\frac{2\RE(\overline{z}G(z))}{1-|z|^2}t+o(t).\]
The condition (\ref{B}) appears as $t$ tends to $0^+$.\\

 We now assume that the condition (\ref{B}) is valid. For $z_0\in\DD$, consider the initial value problem 
\[\frac{dw}{dt}=G(w),\quad w(0)=z_0.\]
Since $G$ is analytic and thus locally Lipschitz, there exist local solutions $w(t)=\varphi_t(z_0)$ by the Cauchy--Peano theorem with values in $\DD$. Let 
\[\rho(z_1,z_2)=\min_{\gamma(0)=z_1;\gamma(1)=z_2}\int_\gamma \frac{2}{1-|z|^2}|dz|.\]
So
\begin{eqnarray*}
\rho(z_0,\varphi_t(z_0))&\leqslant& \int_0^t \frac{2}{1-|\varphi_s(z_0)|^2}\left|\frac{\partial \varphi_s(z_0)}{\partial s}\right|ds \\
& = & \int_0^t \frac{2}{1-|\varphi_s(z_0)|^2}\left|G( \varphi_s(z_0))\right| ds.
\end{eqnarray*}
Write $\displaystyle f:t\mapsto \frac{2}{1-|\varphi_t(z_0)|^2}\left|G( \varphi_t(z_0))\right|$, so that
\begin{eqnarray*}
f'(t)&=& \frac{2}{(1-|\varphi_t(z_0)|^2)^2}\bigg[ \frac{\partial |G(\varphi_t(z_0))|}{\partial t} (1-|\varphi_t(z_0)|^2)\\
& & \qquad\qquad  +2\RE\left( \overline{\varphi_t(z_0)}G(\varphi_t(z_0)) \right)|G(\varphi_t(z_0))|\bigg]\\
&=& \frac{2|G(\varphi_t(z_0))|}{(1-|\varphi_t(z_0)|^2)^2}\Big[ \RE(G'(\varphi_t(z_0))) (1-|\varphi_t(z_0)|^2)\\
& & \quad\qquad\qquad\qquad\qquad+2\RE\left( \overline{\varphi_t(z_0)}G(\varphi_t(z_0)) \right)\Big]\\
&\leqslant &0 \quad \text{by condition (\ref{B}) at $\varphi_t(z_0)$}.
\end{eqnarray*}
We conclude that $f$ is a decreasing function, and thus, for $0\leqslant t_1 <t_2<\eta$,

\[\rho(\varphi_{t_1}(z_0),\varphi_{t_2}(z_0))\leqslant (t_2-t_1) \frac{2|G(\varphi_{t_1}(z_0))|}{1-|\varphi_{t_1}(z_0)|^2}.\]
Therefore, on $[0, \eta)$, $\varphi_t$ remains in a compact subset of $\DD$, so
\[\rho(\varphi_{t_1}(z_0),\varphi_{t_2}(z_0))\leqslant K|t_2-t_1|,\]
where $K$ is a constant independent of $t_1,t_2$ for $0\leqslant t_1 <t_2<\eta$.
Thus, $\varphi_t(z_0)$ converges as $t$ tends to $\eta$. This proves that there exists a solution on $\RR_+$ of the initial value problem. Following \cite{B_P}, $A$ generates a $C_0$-semigroup of composition operators on $H^2(\DD)$.

\endpf

\begin{rem}{\rm
A condition similar to \eqref{B} appears in the paper \cite{irish}, expressed in the language of 
semi-complete vector fields (semiflows); that is, solutions to the Cauchy problem
\begin{eqnarray*}
\frac{du}{dt}+f(u) &=& 0, \\
u(0) &=& x,
\end{eqnarray*}
together with the alternative condition
\[
\re f(z)\overline z \ge \re f(0)\overline z (1-|z|^2), \qquad z \in \DD,
\]
(see also \cite[Prop. 3.5.2]{shoikhet}). Thus, as in \cite{B_P}, they start with
a semigroup of  functions under composition.
}
\end{rem}

\begin{Not}
For each $G(z)=\sum_{n=0}^\infty \alpha_n z^n \in H^2(\beta)$, we write $\tilde G(z) = \alpha_1 + (\alpha_2+\overline{\alpha_0})z+\sum_{n=3}^\infty \alpha_n z^{n-1}$. 
\end{Not}

An easy test using numerical ranges gives the following necessary condition for the generation of a $C_0$ semigroup
of quasicontractions. A more general result (with a more complicated proof) appears in
Proposition~\ref{prop:tinto}.

\begin{prop}
If the operator $A$ generates a $C_0$-semigroup of quasicontractions on $H^2(\beta)$ with $\beta =\left(n^{-\alpha}\right)$ and $\alpha \ge 0$, then
\begin{equation}\label{C}
\mathop{\rm ess}\sup_{z\in\TT} \RE (\tilde G(z))=\mathop{\rm ess}\sup_{z\in \TT} \RE (\overline{z} G(z))\leqslant 0.
\end{equation}
\end{prop}

\beginpf Observing that 
\[\sup_{z\in \TT} \RE (\overline{z} G(z)) =
\sup_{\theta\in\RR} \left\{\RE(\alpha_1) +\RE\left((\overline{\alpha_0} +\alpha_2)e^{i\theta} + \sum_{n=3}^{\infty}\alpha_ne^{i(n-1)\theta}\right)\right\},\]
we can compute the numerical range of  $A$. Let $f$ be an analytic
function defined by $f(z)=\sum_{n=0}^\infty a_nz^n$ with
$\|f\|_{H^2(\beta)}=1$ and $f\in D(A)$. Then we have 
\begin{eqnarray*}
\RE\left(\left\langle G(z)f'(z),f(z)\right\rangle\right)
&=&\RE\left(\left\langle \tilde G(z)zf'(z),f(z)\right\rangle
+\overline{ \alpha_0} \sum_{n=0}^\infty \beta_n\beta_{n+1}a_{n} \overline{a_{n+1}} \right)\\
& =& \RE(\alpha_1)\sum_{n=0}^\infty\beta_{n}^2 n |a_n|^2  \\
 &  & +\RE\left((\alpha_2 +\overline{\alpha_0})\sum_{n=1}^\infty\beta_n\beta_{n+1} n a_n \overline{a_{n+1}} \right)\\
& & + \RE\left(\overline{ \alpha_0} \sum_{n=0}^\infty \beta_n\beta_{n+1}a_{n} \overline{a_{n+1}} \right)\\
& & +  \RE\left(\sum_{k=3}^{\infty}\alpha_k \sum_{n=0}^\infty \beta_n\beta_{n+k-1}n a_n \overline{a_{n+k-1}}\right).
\end{eqnarray*}
Consider 
the polynomial functions
(obviously in $D(A)$) defined by 
\[f_N(z)=c_N\sum_{n=1}^N \frac{\sqrt{6}e^{-in\theta}}{\pi
  n^{1-\alpha}}z^n,\]
where $c_N$ is  a positive real chosen so that $\|f_N\|_{H^2(\beta)}
=1$. It is clear that $c_N$ tends to $1$ as $N$ tends to $\infty$. 
Note now that if (\ref{C}) is not satisfied,  then a suitable choice of $\theta$ makes
$\RE\left( \langle Af_N,f_N\rangle\right)$
tend to $\infty$ as $N$ tends to $\infty$.  It follows that if
(\ref{C}) is not satisfied, then $A$ cannot generate a $C_0$-semigroup, see
e.g. \cite[Chap.II, 3.23]{E_N}.

\endpf

\begin{rem}
It is easy to see that condition (\ref{B}) implies condition (\ref{C}).
\end{rem}

\begin{prop}\label{prop:CimpliesB}
The condition (\ref{C}) implies the condition (\ref{B}).
\end{prop}
\beginpf
Assume $G$ satisfies condition (\ref{C}) and let $H(z)=z\tilde G(z)$. Condition (\ref{C}) and the maximum principle implies that
 $\sup_{z\in\DD}\tilde G(z)\leqslant 0$. Thus by \cite[Theorem 3.3]{B_P} and Corollary~\ref{cor:quasicont}, $f\mapsto Hf'$ generates a $C_0$-semigroup of composition operators on $H^2(\DD)$. Hence, by Theorem \ref{ThmB}, $\forall z\in\DD$,
\[ X:=2\RE(\overline z H(z))+(1-|z|^2)\RE H'(z)\leqslant 0.\]
Now
\begin{eqnarray*}
X &=& 
\RE((1+|z|^2)a_1+2(\overline{a_0}+a_2)z + \sum_{k=3}^\infty a_k z^{k-1}(k-(k-2)|z|^2) \\
&=& \RE(2 a_0 \overline z+(1+|z|^2)a_1+2a_2z + \sum_{k=3}^\infty a_k z^{k-1}(k-(k-2)|z|^2) \\
&=& \RE\left(2\left(a_0 \overline z + \sum_{k=1}^\infty a_k z^{k-1}|z|^2\right)+(1-|z|^2)\sum_{k=1}^\infty ka_k z^{k-1}\right)\\
&=& 2\RE(\overline z G(z))+(1-|z|^2)\RE G'(z),
\end{eqnarray*}
giving condition (\ref{B}).

\endpf

\begin{prop}\label{prop:tinto}
Let $(\beta_n)_n$ be a decreasing sequence of positive reals such that $\liminf_{n\to\infty}|\beta_n|^{1/n}\geq 1$ and 
let $G\in H^2(\beta)$ such that 
\[\mathop{\rm ess}\sup_{w\in \TT}\RE (\overline{w}G(w))>0.\]
Then \[\sup\RE \{\langle Af,f\rangle :f\in D(A),\|f\|_{H^2(\beta)}=1\}=+\infty,\]
where $A$ is defined on $D(A)=\{f\in H^2(\beta):Gf'\in H^2(\beta)\}$ by $Af=Gf'$.  
\end{prop}
Before proceeding to the proof, we state the following technical lemma which explains the hypothesis on monotonicity of $(\beta_n)_n$. 
\begin{lem}\label{lem:decreasing}
Let $(\beta_n)_n$ be a decreasing sequence of positive reals. Then for all positive integer $N$, there exists $\eta=\eta(N)>0$ such that for all $z\in\DD$ with 
$|w|>1-\delta$, we have 
\[\sum_{n=0}^N\frac{|w|^{2n}}{\beta_n^2}<   \sum_{n=N+1}^\infty \frac{|w|^{2n}}{\beta_n^2}.\]
\end{lem}
\beginpf
Since $(1/\beta_n)_n$ is increasing, we have 
\[\sum_{n=0}^N\frac{|w|^{2n}}{\beta_n^2}\leq \frac{1}{\beta_N^2}(1+|w|^2+\cdots +|w|^{2N})=\frac{1}{\beta_N^2}\left( \frac{1-|w|^{2N+2}}{1-|w|^2}\right).\]
On the other hand, we have 
\[\sum_{n=N+1}^\infty\frac{|w|^{2n}}{\beta_n^2}\geq \frac{1}{\beta_{N+1}^2}\sum_{n=N+1}^\infty |w|^{2n}=\frac{|w|^{2N+2}}{\beta_{N+1}^2(1-|w|^2)}\geq 
\frac{|w|^{2N+2}}{\beta_N^2(1-|w|^2)}.\]
Since $1-|w|^{2N+2}<|w|^{2N+2}$ is equivalent to $|w|>(1/2)^{1/(2N+2)}$, for all $w\in\DD$ such that $|w|>\eta(N)$ with $\eta(N)=1-(1/2)^{1/(2N+2)}$, we have 
\[\sum_{n=0}^N\frac{|w|^{2n}}{\beta_n^2}<   \sum_{n=N+1}^\infty
\frac{|w|^{2n}}{\beta_n^2}.\]

\endpf

\noindent\textbf{Proof of Proposition~\ref{prop:tinto}}: 
By hypothesis, there exists  $\delta>0$ and a sequence 
 $(w_k)_k\subset \DD$ such that $|w_k|\to 1$ and $\RE (\overline{w_k}G(w_k))\geq \delta$. 
Moreover the condition  $\liminf_{n\to\infty}|\beta_n|^{1/n}\geq 1$ guarantees that the space $H^2(\beta)$ has reproducing kernels 
$k_w$ for all $w\in\DD$. Now consider the sequence $(\widehat{k_{w_k}})_k$ 
of normalized reproducing kernels associated with $(w_k)_k$,
i.e. $\widehat{k_{w_k}}=\frac{k_{w_k}}{\|k_{w_k}\|_{H^2(\beta)}}$.
First assume that $k_{w_k}\in D(A)$. In this case, 
the remainder of the proof consist in checking that 
\[\lim_{k\to \infty}\RE\left( \langle A\widehat{k_{w_k}},  \widehat{k_{w_k}}\rangle_{H^2(\beta)} \right)=+\infty.\]
Note that 
\[\langle A\widehat{k_{w_k}},  \widehat{k_{w_k}}\rangle_{H^2(\beta)} =\langle G(\widehat{k_{w_k}})', \widehat{k_{w_k}}\rangle_{H^2(\beta)} =
\frac{1}{\|k_{w_k}\|^2_{H^2(\beta)}}G(w_k)k'_{w_k}(w_k),\]
where $k'_{w_k}(z)=\sum_{n\geq 1}\frac{n\overline{w_k}^n}{\beta_n^2}z^{n-1}$.  It follows that 
\[\langle Ak_{w_k},k_{w_k}\rangle_{H^2(\beta)}=\sum_{n\geq 1}\frac{nG(w_k)\overline{w_k}|w_k|^{2(n-1)}}{\beta_n^2},\]
and thus 
\[\langle A\widehat{k_{w_k}},  \widehat{k_{w_k}}\rangle_{H^2(\beta)} =\frac{\overline{w_k}G(w_k)}{|w_k|^2}
\frac{\sum_{n\geq 1}\frac{n|w_k|^{2n}}{\beta_n^2}}{\sum_{n\geq 0}\frac{|w_k|^{2n}}{\beta_n^2}}.\]
Now, for all positive integer $N$, take $\eta(N)$ as in Lemma~\ref{lem:decreasing}, and $k$ sufficiently large so that $|w_k|>1-\eta(N)$. 
Then we have    
\begin{eqnarray*}
\frac{\sum_{n\geq 1}\frac{n|w_k|^{2n}}{\beta_n^2}}{\sum_{n\geq 0}\frac{|w_k|^{2n}}{\beta_n^2}} & = & 
\frac{\sum_{n=0}^N\frac{n|w_k|^{2n}}{\beta_n^2}+ \sum_{n=N+1}^\infty \frac{n|w_k|^{2n}}{\beta_n^2}}{\sum_{n\geq 0}^N
\frac{|w_k|^{2n}}{\beta_n^2} + \sum_{n=N+1}^\infty\frac{|w_k|^{2n}}{\beta_n^2}}\\
 & \geq & \frac{(N+1)\sum_{n=N+1}^\infty \frac{|w_k|^{2n}}{\beta_n^2}}{2\sum_{N+1}^\infty\frac{|w_k|^{2n}}{\beta_n^2}}=\frac{N+1}{2}.
\end{eqnarray*}
Therefore, for $k$ sufficiently large (so that $|w_k|>1-\eta(N)$), we get 
\[\RE\left( \langle A\widehat{k_{w_k}},  \widehat{k_{w_k}}\rangle_{H^2(\beta)} \right)\geq \frac{(N+1)}{2|w_k|^2}\RE(\overline{w_k}G(w_k)).\]
Since $\RE (\overline{w_k}G(w_k))\geq \delta$ and since $|w_k|$ tends to $1$, we get the desired conclusion. 

If $k_{w_k}$ is not in $D(A)$, the conclusion follows from similar
calculation, considering the sequence of polynomials
$(k^M_{w_k})_{M\geq 0}$
defined by 
\[ k^M_{w_k}=\sum_{n=0}^M\frac{\overline{w_k^n}}{\beta_n^2}z^n,\]
which belongs to $D(A)$ and tends to $k_{w_k}$ in $\D$. 

\endpf

We are now ready for the main theorem of this section.

\begin{thm}\label{thm:hardy}
Let $G \in H^2(\DD)$ and $A$ the operator $f \mapsto Gf'$, defined on the
domain $D(A)=\{f \in H^2(\DD): Gf' \in H^2(\DD)\}$ which is dense in $H^2(\DD)$. Then the following conditions are equivalent:\\
(i) $A$ generates a $C_0$-semigroup of composition operators on $H^2(\DD)$;\\
(ii) $2\RE \overline z G(z) + (1-|z|^2) \RE G'(z) \le 0$ for $z \in \DD$;\\
(iii) $A$ generates a quasicontractive $C_0$-semigroup on $H^2(\DD)$;\\
(iv) $\esssup_{z \in \TT} \RE \overline z G(z)  \le 0$.
\end{thm}

\beginpf
The equivalence between (i) and (ii) is   Theorem~\ref{ThmB}. 
The implication (i)$\Rightarrow$(iii) is Corollary~\ref{cor:quasicont}.
The implication (iii)$\Rightarrow$(iv) follows immediately from Proposition~\ref{prop:tinto}  with $\beta_n=1$ for all $n$.
Finally, the implication (iv)$\Rightarrow$(ii) is Proposition~\ref{prop:CimpliesB}.

\endpf

Let $A$ be defined on $D(A):=\{f\in H^2(\DD): Gf'\in H^2(\DD)\}$ by  $A f(z) = G(z) f'(z)$ where $G(z)=\sum_{n=0}^{\infty}\alpha_nz^n$.
An easier condition than Condition (\ref{B}) to test is the following:
\begin{equation}\label{A}
\RE(\alpha_1) +|\overline{\alpha_0} +\alpha_2| + \sum_{n=3}^{\infty}|\alpha_n|\leqslant 0.
\end{equation}
In the sequel we present the link between Condition (\ref{A}) and
Condition (\ref{B}). 
\begin{prop}~%\\
\begin{enumerate}[(i)]
\item Condition (\ref{A}) implies Condition (\ref{B}).
\item If $G\in \CC_2[X]$ (i.e., a polynomial of degree at most 2), then conditions (\ref{B}) and (\ref{A}) are equivalent.
\item There exists a polynomial function of degree $3$ such that condition (\ref{B}) holds and condition (\ref{A}) does not.
\end{enumerate}
\end{prop}

\beginpf
\begin{enumerate}[(i)]
\item The condition (\ref{B}) is equivalent to
\[ (1+|z|^2)\RE(\alpha_1) +2\RE((\overline{\alpha_0}+\alpha_2)z) +  \sum_{n=3}^\infty \RE(\alpha_n((2-n)|z|^2+n)z^{n-1})\leqslant 0.\]
The condition (\ref{A}) is 
\[\RE(\alpha_1) +|\overline{\alpha_0} +\alpha_2| + \sum_{n=3}^{\infty}|\alpha_n|\leqslant 0,\]
that is,
\[(1+|z|^2)\RE(\alpha_1) +(1+|z|^2)|\overline{\alpha_0} +\alpha_2| + \sum_{n=3}^{\infty}(1+|z|^2)|\alpha_n|\leqslant 0.\]
Note that $2\RE((\overline{\alpha_0}+\alpha_2)z)\leqslant 2|\overline{\alpha_0} +\alpha_2||z|\leqslant(1+|z|^2)|\overline{\alpha_0} +\alpha_2|$. On the other hand, the arithmetico-geometric inequality gives, $\forall k\in\NN^*,\ \forall x\in\RR_+$,
\[ \frac{1+x^2+(k-1)x^{k+2}}{k+1} \geqslant \sqrt[k+1]{1\times x^2\times x^{(k+2)(k+1)}}=x^k,\]
i.e. \[ 1+x^2 \geqslant x^k((k+1)-(k-1)x^2).\]
We now observe that 
\[\RE( \alpha_n((2-n)|z|^2+n)z^{n-1}) \leqslant |\alpha_n||z^{n-1}|((2-n)|z|^2+n) \leqslant (1+|z|^2)|\alpha_n|.\]

\item Let $G(z) = \alpha_0+\alpha_1z+\alpha_2z^2$.\\
If condition (\ref{B}) is true, we have in particular 
\[\forall \theta\in\RR,\  2\RE(e^{-i\theta}G(e^{i\theta}))\leqslant 0\quad \text{ i.e.}\quad  \RE(\alpha_1+(\overline{\alpha_0}+\alpha_2)e^{i\theta})\leqslant 0.\]
For $\theta=-\operatorname{arg}(\overline{\alpha_0}+\alpha_2)$, we get \[\RE(\alpha_1)+|\overline{\alpha_0}+\alpha_2|\leqslant 0.\]

\item Take $G(z)=-z+\frac{z^2}{\sqrt{3}}-\frac{z^3}{\sqrt{3}}$. Note
  that $G$ does not satisfy condition (\ref{A}). On the other hand,
  note that
  $G(z)=-z(1-\frac{z}{\sqrt{3}}+\frac{z^2}{\sqrt{3}})=-zF(z)$. In
  \cite{B_P}, it is shown that, if $\RE(F)\geqslant 0$, then
  $G$ generates a $C_0$-semigroup of composition operators on
  $H^2(\beta)$, which is equivalent to condition (\ref{B}). Since
  $\RE(F)$ satisfies the maximum principle, for
  $h(\theta):=\RE(F(e^{i\theta})=1-\frac{1}{\sqrt{3}}\cos(\theta)+\frac{1}{\sqrt{3}}\cos(2\theta)$,
  $F$ maps the disc into the right-half plane if  $h$ is
  nonnegative. For that purpose, note that 
\[ h'(\theta)=\frac{1}{\sqrt{3}}\sin(\theta)(1-4\cos(\theta)).\]
It follows that $h'(\theta)=0 \Leftrightarrow \theta =0 \text{ or } \theta = \pi
\text{ or } \cos(\theta)=\frac 14$. A direct computation gives: $h(0)=1$, $h(1)=1+\frac{2}{\sqrt{3}}$ and if $\cos(\theta)=\frac 14$, $h(\theta)=1-\frac{9}{8\sqrt{3}}>0$. Therefore $G$ satisfies condition (\ref{B}).
\end{enumerate}

\endpf

\section{Quasicontractive semigroups on the Dirichlet space}\label{sec:4}

Recall that  the Dirichlet space norm is defined by
\begin{equation}\label{eq:sep29}
\|f\|^2_\D = |a_0|^2 + \sum_{k=1}^\infty k|a_k|^2 = |f(0)|^2 + \int_\DD |f'(z)|^2 \, dA(z),
\end{equation}
for $f(z)=\sum_{k=0}^\infty a_k z^k$,
and it is induced by an inner product that may be written, at least formally, as
\[
\langle f,g \rangle_\D = \langle f,zg' \rangle_{H^2(\DD)} + f(0)\overline{g(0)}.
\]

\begin{prop}\label{prop:doc5}
For $G \in \D$ and $A$ the operator $f \mapsto Gf'$, defined on the
domain $D(A)=\{f \in \D: Gf' \in \D\}$, which in dense in $\D$, the following two conditions are equivalent:\\
(i) $\esssup_{z \in \TT} \RE \overline z G(z)  \le 0$;\\
(ii) $\sup\{ \RE \langle Af,f \rangle_\D : f \in D(A), \|f\|_\D=1 \} < \infty$.
\end{prop}

\beginpf
Suppose that $\esssup_{z \in \TT} \RE \overline z G(z)  \le 0$. 
%Note that $\RE \overline zG(z) = \RE \widetilde G(z)$
%for $z \in \TT$,
%where, for $G(z)=\sum_{k=0}^\infty a_k z^k$ we define
%\[
%\widetilde G(z) = a_1+(a_2+\overline{a_0})z+a_3z^2+a_4z^3+\ldots.
%\]
Then
\begin{eqnarray*}
\RE \langle Af,f\rangle_\D  &=& 
\RE \langle Gf',zf'\rangle_{H^2(\DD)} + \RE\left( G(0)f'(0)\overline{f(0)} \right)\\
&=&
\RE \left( \frac{1}{2\pi}\int_0^{2\pi} G(z) \overline z |f'(z)|^2 \, d\theta\right) + \RE\left( G(0)f'(0)\overline{f(0)}\right),
\end{eqnarray*}
 with $ z=e^{i\theta}$,
and the supremum of this quantity over $\|f\|_\D=1$, $f\in D(A)$ is clearly finite.\\

Conversely, suppose that $\esssup \RE \overline z G(z) > 0$. 
By considering  an (outer) function $u$ with $|u|=1$ on a set of positive measure where $\RE \overline zG(z) > \delta > 0$ and
$|u|=1/2$ on its complement (see Thm. 4.6 in \cite{Garnett})
we see that $\liminf_{n \to \infty}\RE \langle Gu^n,zu^n \rangle_{H^2(\DD)} > 0$.
It now follows that there is a function $f \in D(A)$ with $\langle f,f \rangle_\D=1$
and $\RE \langle Gf',zf' \rangle_{H^2(\DD)} > 0$.

Now define a sequence $(f_k)_k$ in $\D$ by setting $f_k'=z^k f'$ and $f_k(0)=0$.

So if $f(z)=\sum_{n=0}^\infty a_n z^n$, then 
\[
f_k(z)= \sum_{n=1}^\infty a_n \frac{n}{n+k} z^{n+k},
\]
and hence
\begin{eqnarray*}
\langle f_k , f_k \rangle_\D &=&  \sum_{n=1}^\infty |a_n|^2 \left( \frac{n}{n+k} \right)^2 (n+k) \\
&=& 
 \sum_{n=1}^\infty |a_n|^2 \frac{n^2}{n+k} \le \|f\|^2_\D,
\end{eqnarray*}
and thus this tends to zero by the dominated convergence theorem.

Now 
\[
\RE \langle Af_k,f_k \rangle_\D = \RE \langle \overline z G f'_k,f'_k \rangle_{H^2(\DD)} =  \RE    \langle   G f',f' \rangle_{H^2(\DD)} 
\]

On normalizing the functions $f_k$ we see that
\[
\sup\{ \RE \langle Af,f \rangle_\D : f \in \D, \|f\|_\D=1 \} = \infty.
\]

\endpf

\begin{prop}\label{prop:univalent}
Let $G \in\D$ and $A$ the operator $f \mapsto Gf'$, defined on the
domain $D(A)=\{f \in \D: Gf' \in \D\}$, which is dense in $\D$. If $A$ generates a $C_0$-semigroup of composition operators on $\D$, then
this semigroup is quasicontractive.
\end{prop}

\beginpf
Given a semigroup $(C_{\varphi_t})_{t\geq 0}$ acting on $\D$, we must
show that   $\|C_{\varphi_t} f\|_\D = \|f\|_\D (1+ O(t))$ for small $t>0$.
First, since $\phi_t$ is injective, we have the well-known inequality
\begin{eqnarray*}
\int_\DD |(f \circ \phi_t)'(z)|^2 \, dA(z) &=& \int_\DD |(f' \circ \phi_t(z))|^2 |\phi'_t(z)|^2 \, dA(z) \\ 
&= & \int_{\phi(\DD)} |f'(w)|^2 \, dA(w) \le \int_\DD |f'(w)|^2 \, dA(w),
\end{eqnarray*}
taking $w=\phi(z)$. Therefore the composition operator $C_{\phi_t}$ is
bounded on $\D$. Moreover,  by \cite[Thm.~2]{martin05}, 
\[\|C_{\varphi_t}\|\leq \sqrt{   \frac{L+2+\sqrt{L(4+L)}}{2}} 
,\]
where $L=\log\left( \frac{1}{1-|\varphi_t(0)|^2} \right)$. This upper
bound is sharp since it is an equality whenever $\DD\setminus
\varphi_t(\DD)$ is of Lebesgue area measure equal to $0$. 
 
Siskakis  \cite{sisk96} proved that, as in the case of the Hardy space, $A$ is of the form $A(f)=G(z)f'(z)$,
where $G$ is an holomorphic function on $\DD$ and
$\varphi_t(z)=z+G(z)t+o(t)$.   It follows that, for $t\to 0$, 
 \[\|C_{\varphi_t}\|\leq 1+O(t),\]
 since  $L=O(t^2).$
Therefore, there exists $w>0$ such that $\|C_{\varphi_t}\|\leq e^{wt}$
for all $t\geq 0$, and thus $(C_{\varphi_t})_{t\geq 0}$ is then a
quasicontractive $C_0$-semigroup.

\endpf

\begin{thm}\label{thm:dir}
Let $G \in\D$ and $A$ the operator $f \mapsto Gf'$, defined on the
domain $D(A)=\{f \in \D: Gf' \in \D\}$, which is dense in $\D$. Then the following conditions are equivalent:\\
(i) $A$ (extended to its natural domain in $H^2(\DD)$) generates a $C_0$-semigroup of composition operators on $H^2(\DD)$;\\
(ii) $2\RE \overline z G(z) + (1-|z|^2) \RE G'(z) \le 0$ for $z \in \DD$;\\
(iii) $A$ (extended to its natural domain in $H^2(\DD)$) generates a quasicontractive $C_0$-semigroup on $H^2(\DD)$;\\
(iv) $\esssup_{z \in \TT} \RE \overline z G(z)  \le 0$;\\
(v) $\sup\{ \RE \langle Af,f \rangle_\D : f \in D(A), \|f\|_\D=1 \} < \infty$;\\
(vi) $A$ generates a quasicontractive $C_0$-semigroup on $\D$;\\
(vii) $A$ generates a $C_0$-semigroup of composition operators on $\D$.
\end{thm}

\beginpf
Conditions (i)--(iv) have already been shown to be equivalent in Theorem~\ref{thm:hardy}.
The equivalence of conditions (iv) and (v) is shown in Proposition~\ref{prop:doc5}.
For (i)$\Rightarrow$(vii) is detailed in \cite{sisk96}.
The fact that (vii)$\Rightarrow$(vi) is given in Proposition~\ref{prop:univalent}.
Finally, (vi)$\Rightarrow$(v) by Lumer--Phillips result (see 
Corollary~\ref{cor:lm} ).

\endpf

\section{Comments}
In \cite{B_P}, as well as in Condition~\ref{A}, the description of the
generator of a $C_0$-semigroup of composition  operators 
 relies on analytic
functions $F$ or $\tilde{G}$ which map $\DD$ into the right or left half-plane.    
For that purpose, let us recall the  Carath\'eodory--Toeplitz theorem \cite{B_C,S}.  

 \begin{thm}[Carath\'eodory--Toeplitz]
 Let $f(z)=\sum_{n=0}^\infty \mu_nz^n$ and consider for $k\geqslant 1$ the matrices $M_k=(m_{i,j})_{1\leqslant i,j \leqslant k}$ where $m_{i,j}=\mu_{j-i}$ if $i\leqslant j$ and $m_{i,j}=0$ otherwise. Then $f$ maps the disc to the right half plane if and only if the Hermitian matrix $N_k=M_k+\overline{M_k'}$ is nonnegative definite for all $k\geqslant 1$.
 \end{thm}
 This theorem has to be considered with the Sylvester Criterion.
 \begin{thm}[Sylvester Criterion]
 Let $A=\big(a_{ij}\big)_{1\leqslant i,j\leqslant n}$ be Hermitian. Then $A$ is positive definite if and only if the $n$ matrices $A_p=\big(a_{ij}\big)_{1\leqslant i,j\leqslant p}$ with $1\leqslant p \leqslant n$ have positive determinant.
 \end{thm}
Here is an example where we can use those tools. 
 \begin{ex}
 Let $G(z)=a_0+a_1z+a_2z^2 \in\CC_2[X]$, thanks to condition (\ref{C}), we have that $G$ generates a $C_0$-semigroup of analytic functions on $\DD$ if and only if $\sup(\tilde{G}(\T))\leqslant 0$. Besides, $\sup(\tilde{G}(\TT))< 0$ if and only if
 \[\det \begin{pmatrix} -\RE(a_1) \end{pmatrix} >0,\quad \det\begin{pmatrix}
 -\RE(a_1) & -(\overline{a_0}+a_2) \\ -(\overline{a_2}+a_0)&-\RE(a_1)
 \end{pmatrix} >0\]
 i.e. 
 \[ \left\{ \begin{array}{l}
 \RE(a_1)<0\\
 \RE(a_1)^2-|\overline{a_0}+a_2|^2=(\RE(a_1)-|\overline{a_0}+a_2|)(\RE(a_1)+|\overline{a_0}+a_2|)>0
 \end{array}
 \right.\]
 i.e.
 \[ \RE(a_1)+|\overline{a_0}+a_2|<0.\]
 We have recovered condition (\ref{A}).
 \end{ex}

% %\subsection{Examples}

One may wonder if the quasicontractive $C_0$-semigroup whose generator
is given by $A$ can be determined on $H^2(\DD)$ or $\D$. We know that
it is a semigroup of  composition operators $C_{\varphi_t}$, with
\[
\frac{\partial \varphi_t(z)}{\partial t}=G(\varphi_t(z)).
\]
This is an important and not so easy issue, which can be answered in
some particular cases, as follows. Those examples are slight
generalizations of the ones presented in \cite{AS}. 

\begin{enumerate}[(i)]
\item If $G(z)=az+b$ with $a\neq 0$, we have 
 \[\varphi_t(z)= e^{at}z+\frac{b}{a}(e^{at}-1).\]
Furthermore, the Denjoy--Wolff point $\alpha$ of this holomorphic semigroup is $\alpha = -\frac ba\in\DD$.
\item When \emph{$G$ is a polynomial of degree $2$}, defined by $G(z)=c(z-a)(z-b)$:
\begin{itemize}
\item If $a\neq b$, we get
\[
\varphi_t(z) = \frac{z(ae^{bct}-be^{act})+ab(e^{act}-e^{bct})}{z(e^{bct}-e^{act})+(ae^{act}-be^{bct})},
\]
whose Denjoy--Wolff point is $\alpha = a\in\DD$ if $\RE a<\RE b$ and
$\alpha = b\in\DD$ if $\RE a>\RE b$.
 In the case where $\RE a =\RE b$ it happens that $\varphi_{t_n}=Id$ for $t_n=\frac{2\pi n}{\IM a-\IM b}$, so $\varphi_t$ is an automorphism.
\item If $a=b$, we find another expression for $\varphi_t$:
\[
\varphi_t(z) = \frac{z(1-act)+a^2ct}{-zct+(1+act)},
\]
whose Denjoy--Wolff point is $\alpha=a$.
\end{itemize}
\item As \emph{$G$ is polynomial of higher degree}, we usually do not have explicit expression of the semigroup $(\varphi_t)$. Yet, some cases can be found:
\begin{itemize}
\item If $G(z)=c(z-a)^n$  then $\forall t\in\RR_+,\forall z\in\DD$, 
\[\varphi_t(z)=a+\frac{z-a}{(1-nct(z-a)^{n-1})^{\frac{1}{n-1}}}.\]
Note that, if $c=1$, the only possible case is when $a=1$.
\item
If  $G(z)=cz(z^n-a)$  then $\forall t\in\RR_+,\forall z\in\DD$, 
\[\varphi_t(z)=\frac{ze^{-ct}}{\left(1-z^n\left(\frac{1-e^{-nct}}{a}\right)\right)^{\frac{1}{n}}}.\]
\end{itemize}
\end{enumerate}

\subsection*{Acknowledgement}
The authors are grateful to the referee for drawing their
attention to the book of Shoikhet \cite{shoikhet}, and hence
indirectly to \cite{irish}.

\end{document}